# On some further properties and application of Weibull-$R$ family of distributions

Indranil Ghosh[1], Saralees Nadarajah[2]
[1]University of North Carolina, Wilmington, USA
[2] University of Manchester, UK

## Abstract

In this paper, we provide some new results for the Weibull-R family of distributions (Alzaghal, Ghosh and Alzaatreh (2016)). We derive some new structural properties of the Weibull-R family of distributions. We provide various characterizations of the family via conditional moments, some functions of order statistics and via record values.

## 1 Introduction

Various well known univariate distributions have been extensively used over the past few decades for modeling data arising from different spheres such as engineering, actuarial, environmental and medical sciences, biological studies, demography, economics, finance and insurance. However, in many of these applied areas in particular lifetime analysis, finance and insurance, there is a growing demand for extended forms of these distributions. As a consequence, several methods for generating new families of distributions have been studied in the literature.

The Weibull distribution is a well-known distribution due to its extensive use to model various types of data. This distribution has been widely used in survival and reliability analyses. This distribution has quite a bit of flexibility for analyzing skewed data. It allows for increasing and decreasing hazard rate functions (hrfs), depending on the shape parameters, which gives an extra edge over the exponential distribution that has only constant hrf. Since the 1970s, many extensions of the Weibull distribution have been proposed to enhance its capability to fit diverse lifetime data and Murthy et al. (2004) proposed a scheme to classify these distributions.

Unfortunately Weibull distribution has certain drawbacks. Bain (1978) pointed out that the maximum likelihood estimators of the Weibull parameters may not behave properly for all parameters even when the location parameter is zero. When the shape parameter is greater than one, the hrf of the Weibull distribution increases from zero to infinity, which may not be appropriate in some scenarios. Also, the Weibull family does not enjoy the likelihood ratio ordering, as a consequence, there does not exist a uniformly most powerful test for testing one-sided hypotheses on the shape parameter. The sum of independently and identically distributed Weibull random variables is not difficult to obtain, one may use the convolution and/or characteristic approach to find the distribution of the sum which will involve some special functions. Mudholkar et al. (1995) proposed a three parameter (one scale and two shape) distribution, the exponentiated Weibull (see also Mudholkar and Srivastava (1993)) distribution. They indicated, based on certain data sets, the exponentiated Weibull distribution provides better fits than the two parameter Weibull distribution. Gupta and Kundu (1997) considered a special case of the exponentiated Weibull distribution assuming location parameter to be zero. They compared its performance with the two parameter gamma family of distributions and the two parameter Weibull family, through data analysis and computer simulations.

In quest for a greater applicability of the Weibull distribution many researchers have considered various types of generalizations. These generalizations include broad family of univariate distributions generated from the Weibull distribution introduced by Gurvich et al. (1997), the generalized Weibull distribution due to Mudholkar and Kollia (1994), and the beta-Weibull distribution due to Lee et al. (2007). The log-Weibull distribution has been studied in detail by many authors. For example, White (1969) studied the moments of log-Weibull order statistics, while Huillet and Raynaud (1999) studied their application in earthquake magnitude data. Ortega et al. (2013) discussed usefulness of the log-Weibull regression model to predict recurrence of prostate cancer.

Generalized Weibull distributions can be constructed in many ways, as detailed in Lai et al. (2011) and references therein. Members of this family usually contain the standard Weibull distribution as a special case.

Nadarajah and Kotz (2005) defined a class of extended Weibull (EW) distributions with cumulative distribution function (cdf) given by $G_{\alpha,\boldsymbol{\tau}}(t) = 1 - \exp\{-\alpha H(t)\}$, where $\alpha > 0$ and $H(t)$ is a monotonically increasing function of $t$ with the only limitation $H(t) \geq 0$ and $\boldsymbol{\tau}$ represents a vector of unknown parameters in $H(t)$. If $H(t)$ is a power law function, the above equation reduces to the traditional Weibull distribution. In this paper, we consider a new generalization of any absolutely continuous (R) distribution, using Weibull as a baseline distribution, called the Weibull-$R$ family of distributions, following the technique of Alzaatreh et al. (2013b). It is to be noted here that parallel development for any discrete distribution with the baseline distribution as Weibull can also be developed. It is defined as follows:

Let $T \in (a, b)$, $R$ and $Y \in (c, d)$ be random variables with cdfs $F_T(x) = P(T \leq x)$, $F_R(x) = P(R \leq x)$ and $F_Y(x) = P(Y \leq x)$ for $-\infty \leq a < b \leq \infty$ and $-\infty \leq c < d \leq \infty$. Here, $R$ can be a continuous or a discrete random variable. Let $Q_T(p)$, $Q_R(p)$ and $Q_Y(p)$ denote the corresponding quantile functions, where the quantile function of a random variable $Z$ is defined as $Q_Z(p) = \inf\{z : F_Z(z) \geq p\}$, $0 < p < 1$. If the probability density functions (pdfs) of $T$, $R$ and $Y$ exist, we denote them by $f_T(x)$, $f_R(x)$ and $f_Y(x)$, respectively. We define a random variable $X$ as having the cdf

$$F_X(x) = \int_a^{Q_Y(F_R(x))} f_T(t)dt = F_T\left(Q_Y\left(F_R(x)\right)\right) \tag{1}$$

for $-\infty < x < \infty$. Alzaatreh et al. (2013b) referred to the distributions in (1) as the $T$-$RY$ family of distributions. The pdf and hazard rate funtion (hrf) of $X$ can be derived as

$$f_X(x) = f_R(x) \cdot \frac{f_T\left(Q_Y\left(F_R(x)\right)\right)}{f_Y\left(Q_Y\left(F_R(x)\right)\right)} \tag{2}$$

and

$$h_X(x) = h_R(x) \cdot \frac{h_T\left(Q_Y\left(F_R(x)\right)\right)}{h_Y\left(Q_Y\left(F_R(x)\right)\right)}.$$

We can rewrite (1) and (2) as

$$F_X(x) = F_T\left(-\log\left(1 - F_R(x)\right)\right) = F_T\left(H_R(x)\right)$$

and

$$f_X(x) = \frac{f_R(x)}{1 - F_R(x)} f_T\left(-\log\left(1 - F_R(x)\right)\right) = h_R(x) f_T\left(H_R(x)\right), \quad (3)$$

where $h_R(x)$ and $H_R(x) = -\log\left(1 - F_R(x)\right)$, where, $h_R(x)$ is the hazard rate function for the random variable $R$, and $H_R(x)$ being the survival function of $R$. The cdf of a random variable $X$ can take this form only if a random varibale $Y$ is unit exponentially distributed, i.e., $Y \sim Exp(1)$.

If $T$ is a Weibull random variable with parameters $c$ and $\gamma$, (3) is the pdf of the Weibull-$R$ distribution:

$$f_X(x) = \frac{c}{\gamma} \frac{f_R(x)}{1 - F_R(x)} \left[\frac{-\log\left(1 - F_R(x)\right)}{\gamma}\right]^{c-1} \exp\left\{-\left[\frac{-\log\left(1 - F_R(x)\right)}{\gamma}\right]^c\right\} \quad (4)$$

for $c > 0$ and $\gamma > 0$.

The cdf corresponding to (4) is

$$F_X(x) = 1 - \exp\left\{-\left[\frac{-\log\left(1 - F_R(x)\right)}{\gamma}\right]^c\right\}. \quad (5)$$

Note that if $R$ is a Weibull random variable then (4) is the pdf of a generalized gamma distribution. Hence, the Weibull-$R$ family is a broad class of distributions as compared to gamma-generalized distributions. For details on the construction see Alzaatreh & Ghosh (2015).

A particular case of the Weibull-$R$ family that we shall study in some detail later is the Weibull-Lomax distribution (WLD), see Sections 3, 5 and 6. Possible shapes of the pdf and the hrf of $X$ for the WLD are shown in Figures 1 and 2. Figure 1 shows that the pdf can be monotonically decreasing or unimodal. Both the left and right tails of the pdf decay to zero slowly, i.e., both tails are heavy. Processes commonly encountered in practice have heavy tails. So, the tails of the WLD are realistic. The right tail of the Weibull distribution decays to zero exponentially, which is not so realistic. Figure 2 shows that the hrf can be monotonically decreasing, monotonically increasing or upside down bathtub shaped. The Weibull distribution cannot exhibit upside down bathtub shaped hrfs. Reliability and survival analysis often encounter upside down bathtub hazard rates. Examples can be found in redundancy allocations in systems (Singh and Misra, 1994) and mortality modeling (Silva et al., 2010).

The contents of this paper are organized as follows. Mathematical properties of the Weibull-$R$ distribution (including characterizations, quantiles, shape properties, entropy measures, moments and reliability parameter) are derived in Sections 2 and 4. Some particular cases of the Weibull-$R$ distribution are studied in Section 3. Finally, Section 5 concludes the paper.

## 2 Characterizations of the Weibull-$R$ family

It is a natural requirement that in designing a stochastic model for a particular modeling problem, an investigator will be vitally interested to know if their model fits the requirements of a specific underlying probability distribution. To this end, the investigator will rely on characterizations of the selected distribution. Generally speaking, the problem of characterizing a distribution is an important problem in various fields and has recently attracted the attention of many researchers. Consequently, various characterizations have been reported in the literature. These characterizations have been established in many different directions. Here, we present characterizations of the newly introduced Weibull-$R$ family of distributions. These characterizations are based on record values. We would like to remark here that other possible ways of characterization of this Weibull-$R$ family might be possible, but, in this present article, we report the most interesting one.

## 2.1 Characterizations of the Weibull-$R$ family via records

Here, we present characterizations of the newly introduced Weibull-$R$ family of distributions via record values.

Let $X_{U(m)}$ and $X_{U(n)}$ for $m < n$ denote the upper record values from a given family specified by pdf $f_X$ and cdf $F_X$. The joint pdf of $X_{U(m)}$ and $X_{U(n)}$ is (Ahsanullah, 1995)

$$f_{X_{U(m)},X_{U(n)}}(x,y) = \frac{[\log(1-F_X(x)) - \log(1-F_X(y))]^{n-m-1}}{\Gamma(m)\Gamma(n-m)}[-\log(1-F_X(x))]^{m-1}\frac{f_X(x)f_X(y)}{1-F_X(x)}, \quad (6)$$

where $-\infty < x < y < \infty$ and $1 \leq m < n$.

**Theorem 1.** If $X \sim$ Weibull-$R(c,\gamma)$, then the pdf of $X_{U(m)}$ is

$$\begin{aligned} f_{X_{U(m)}}(x) &= \frac{f_X(x)}{[1-F_X(x)]}\Gamma(m)\Gamma(n-m)\sum_{j=0}^{n-m-1}(-1)^j\binom{n-m-1}{j}\gamma^{-cj} \\ &\quad \cdot \left[\frac{-\log(1-F_R(x))}{\gamma}\right]^{c(n-2-j)}\Gamma(1+jc, -\log(1-F_R(x))) \end{aligned} \quad (7)$$

for $-\infty < x < \infty$, where $\Gamma(a,x) = \int_x^\infty t^{a-1}\exp(-t)dt$.

*Proof.* From (4), we have $1 - F_X(x) = \exp\left\{-\left[\frac{-\log(1-F_R(x))}{\gamma}\right]^c\right\}$. So,

$$\left[\log\frac{1-F_X(x)}{1-F_X(y)}\right]^{n-m-1} = \left\{\left[\frac{-\log(1-F_R(y))}{\gamma}\right]^c - \left[\frac{-\log(1-F_R(x))}{\gamma}\right]^c\right\}^{n-m-1}.$$

By (6), we can write

$$\begin{aligned} f_{X_{U(m)},X_{U(n)}}(x,y) &= \frac{1}{\Gamma(m)\Gamma(n-m)}\frac{f_X(x)f_X(y)}{1-F_X(x)}\left[\frac{-\log(1-F_R(x))}{\gamma}\right]^{c(m-1)} \\ &\quad \cdot \left\{\left[\frac{-\log(1-F_R(y))}{\gamma}\right]^c - \left[\frac{-\log(1-F_R(x))}{\gamma}\right]^c\right\}^{n-m-1} \end{aligned}$$

for $-\infty < x < y < \infty$. Therefore, the marginal pdf of $X_{U(m)}$ is

$$\begin{aligned} f_{X_{U(m)}}(x) &= \frac{1}{\Gamma(m)\Gamma(n-m)}\frac{f_X(x)}{1-F_X(x)}\left[\frac{-\log(1-F_R(x))}{\gamma}\right]^{c(m-1)} \\ &\quad \cdot \int_x^\infty f_X(y)\left\{\left[\frac{-\log(1-F_R(y))}{\gamma}\right]^c - \left[\frac{-\log(1-F_R(x))}{\gamma}\right]^c\right\}^{n-m-1}dy \\ &= \frac{1}{\Gamma(m)\Gamma(n-m)}\frac{f_X(x)}{1-F_X(x)}\left[\frac{-\log(1-F_R(x))}{\gamma}\right]^{c(m-1)}I_1, \end{aligned} \quad (8)$$

where

$$
\begin{aligned}
I_1 &= \int_x^\infty f_X(y) \left\{ \left[ \frac{-\log(1-F_R(y))}{\gamma} \right]^c - \left[ \frac{-\log(1-F_R(x))}{\gamma} \right]^c \right\}^{n-m-1} dy \\
&= \sum_{j=0}^{n-m-1} (-1)^j \left[ \frac{-\log(1-F_R(x))}{\gamma} \right]^{c(n-m-1-j)} \binom{n-m-1}{j} \\
&\quad \cdot \int_x^\infty f_X(y) \left[ \frac{-\log(1-F_R(y))}{\gamma} \right]^{jc} dy \\
&= \sum_{j=0}^{n-m-1} (-1)^j \gamma^{-cj} \left[ \frac{-\log(1-F_R(x))}{\gamma} \right]^{c(n-m-1-j)} \binom{n-m-1}{j} \\
&\quad \cdot \Gamma(1+cj, -\log(1-F_R(x))). \tag{9}
\end{aligned}
$$

The result follows by substituting (9) in (8). □

Theorem 1 can be useful for estimation based on record values. There are many situations in which only records are observed. Ultimate examples of such situations can be found from the website for Guinness World Records, see http: // www. guinnessworldrecords. com/ . Another example is the situation of testing the breaking strength of wooden beams as described in Glick (1978).

## 3 Some examples of the Weibull-$R$ family

- For $R$ a Pareto random variable with the pdf $f_R(x) = \frac{k\theta^k}{x^{k+1}}$, $x > \theta$, $k > 0$, we have the Weibull-Pareto distribution (WPD) with the pdf, cdf and the hrf given by

$$ f_X(x) = \frac{\beta c}{x} \left[ \beta \log \left( \frac{x}{\theta} \right) \right]^{c-1} \exp \left\{ - \left[ \beta \log \left( \frac{x}{\theta} \right) \right]^c \right\}, $$

$$ F_X(x) = 1 - \exp \left\{ - \left[ \beta \log \left( \frac{x}{\theta} \right) \right]^c \right\} $$

and

$$ h_X(x) = \frac{\beta c}{x} \left[ \beta \log \left( \frac{x}{\theta} \right) \right]^{c-1}, $$

respectively, for $x > \theta$, $c > 0$, $\theta > 0$, and $\beta = k/\gamma$. This family has been studied by Alzaatreh et al. (2013a).

- For $R$ a Lomax random variable with the pdf $f_R(x) = \frac{k}{\theta} \left( 1 + \frac{x}{\theta} \right)^{-k-1}$, $k > 0$ $x > 0$, we have the WLD with the pdf, cdf and the hrf given by

$$ f_X(x) = \frac{\beta c}{x+\theta} \left[ \beta \log \left( 1 + \frac{x}{\theta} \right) \right]^{c-1} \exp \left\{ - \left[ \beta \log \left( 1 + \frac{x}{\theta} \right) \right]^c \right\}, \tag{10} $$

$$ F_X(x) = 1 - \exp \left\{ - \left[ \beta \log \left( 1 + \frac{x}{\theta} \right) \right]^c \right\} \tag{11} $$

and

$$h_X(x) = \frac{\beta c}{x + \theta} \left[ \beta \log \left( 1 + \frac{x}{\theta} \right) \right]^{c-1},$$

respectively, for $x > 0$, $c > 0$, $\theta > 0$ and $\beta = k/\gamma$. Possible shapes of $f_X(x)$ and $h_X(x)$ are shown in Figures 1 and 2.

Note that the WLD is a shifted version of the WPD. When $c = 1$, the WLD reduces to the Lomax distribution with parameters $\beta$ and $\theta$.

- For $R$ a Cauchy random variable with the pdf $f_R(x) = \frac{1}{\pi \left[ 1 + \left( \frac{x}{\delta} \right)^2 \right]}$, $-\infty < x < \infty$, we have the Weibull-Cauchy distribution with the pdf and cdf given by

$$f_X(x) = \frac{2c}{\gamma \left[ 1 + \left( \frac{x}{\delta} \right)^2 \right] \left( \pi - 2 \arctan \frac{x}{\delta} \right)} \left\{ -\frac{\log \left[ \frac{1}{2} - \frac{1}{\pi} \arctan \left( \frac{x}{\delta} \right) \right]}{\gamma} \right\}^{c-1}$$

$$\cdot \exp \left\{ - \left[ -\frac{\log \left( \frac{1}{2} - \frac{1}{\pi} \arctan \left( \frac{x}{\delta} \right) \right)}{\gamma} \right]^c \right\}$$

and

$$F_X(x) = 1 - \exp \left\{ - \left[ -\frac{\log \left( \frac{1}{2} - \frac{1}{\pi} \arctan \left( \frac{x}{\delta} \right) \right)}{\gamma} \right]^c \right\},$$

respectively, for $-\infty < x < \infty$, $\delta > 0$, $\gamma > 0$ and $c > 0$.

- For $R$ a normal random variable with the pdf $f_R(x) = \frac{1}{\sqrt{2\pi}\sigma} \exp \left[ -\frac{1}{2} \left( \frac{x-\mu}{\sigma} \right)^2 \right]$, $-\infty < x < \infty$, we have the Weibull-normal distribution with the pdf and cdf given by

$$f_X(x) = \frac{c}{\gamma \left[ 1 - \Phi \left( \frac{x-\mu}{\sigma} \right) \right]} \frac{1}{\sqrt{2\pi}\sigma} \exp \left[ -\frac{1}{2} \left( \frac{x-\mu}{\sigma} \right)^2 \right]$$

$$\cdot \exp \left\{ - \left[ \frac{-\log \left( 1 - \Phi \left( \frac{x-\mu}{\sigma} \right) \right)}{\gamma} \right]^c \right\} \left[ \frac{-\log \left( 1 - \Phi \left( \frac{x-\mu}{\sigma} \right) \right)}{\gamma} \right]^{c-1}$$

and

$$F_X(x) = 1 - \exp \left\{ - \left[ \frac{-\log \left( 1 - \Phi \left( \frac{x-\mu}{\sigma} \right) \right)}{\gamma} \right]^c \right\},$$

respectively, for $-\infty < x < \infty$, $-\infty < \mu < \infty$, $\sigma > 0$, $\gamma > 0$ and $c > 0$, where $\Phi(\cdot)$ denotes the standard normal cdf.

Pdf graphs will be here.

The WLD, the Weibull-Cauchy distribution and the Weibull-normal distribution are new and do not appear to have been studied by others. The WLD is different from the Weibull-Lomax distribution studied by Tahir et al. (2015), compare (10) with equation (2.1) in Tahir et al. (2015).

## 4 Properties of the Weibull-$R$ family

The corresponding hrf and quantile function for any $u \in (0,1)$ are

$$\begin{aligned}h_X(x) &= \frac{c}{\gamma}\frac{f_R(x)}{1-F_R(x)}\left[\frac{-\log(1-F_R(x))}{\gamma}\right]^{c-1}\\ &= \frac{\partial}{\partial x}\left\{\left[\frac{-\log(1-F_R(x))}{\gamma}\right]^c\right\}\end{aligned}$$

and

$$Q(u) = F_X^{-1}\left(1 - \exp\left\{-\gamma\left[-\log(1-u)\right]^{1/c}\right\}\right).$$

Since these are in closed form, the Weibull-$R$ family can be applied to model censored data also. One can also obtain a closed form expression for the cumulative hrf.

### 4.1 Shape

Lemma 1 below gives the limiting behaviors of the Weibull-$R$ pdf and its hrf. Its proof is obvious.

**Lemma 1.** We have

$$f_X(x) \sim \frac{c}{\gamma^c} f_R(x) F_R^{c-1}(x) \exp\left[-\gamma^{-c} F_R^c(x)\right]$$

and

$$h_X(x) \sim \frac{c}{\gamma^c} f_R(x) F_R^{c-1}(x)$$

as $x \to -\infty$.

**Lemma 2.** The mode of the pdf of the Weibull-$R$ family is the root of

$$\frac{f'_R(x)}{f_R(x)} + \frac{f_R(x)}{1-F_R(x)} - \frac{(c-1)f_R(x)}{[1-F_R(x)]\log[1-F_R(x)]} - c\gamma^{-c}\frac{\{-\log[1-F_R(x)]\}^{c-1}f_R(x)}{1-F_R(x)} = 0,$$

where $f'_R(x) = df_R(x)/dx$.

The proof of Lemma 2 is straightforward.

Analytical solutions to the mode do not appear possible, even for the four examples presented in Section 3. The mode should be computed numerically, for example, using uniroot in the R software.

### 4.2 Entropy measures

Entropies of a random variable $X$ are measures of variation of uncertainty. Entropies have been used in several applications in science, engineering and economics. The Shannon entropy (Shannon, 1951) of a random variable $X$ say with pdf $f_X$ is defined by $E\left[-\log f_X(X)\right]$.

**Lemma 3.** If $X$ is a Weibull-$R$ random variable then its Shannon entropy is

$$\begin{aligned}\eta_X &= -\log c + \log\gamma - (c-1)E\left\{\frac{\log\left([-\log(1-F_R(X))]\right)}{\gamma}\right\} - E\left[\log\frac{f_R(X)}{1-F_R(X)}\right]\\ &\quad + E\left\{\left[\frac{-\log(1-F_R(X))}{\gamma}\right]^c\right\}.\end{aligned}$$

Lemma 3 follows immediately from (4).

The problem of testing whether some given observations can be considered as coming from one of two probability distributions is an old problem in statistics. Consider a random sample $X_1, X_2, \ldots, X_n$ of size $n$ from a Weibull-$R$ family. The objective is to identify a specific Weibull-$R$ distribution in (4) that is most appropriate to describe the data $X_1, X_2, \ldots, X_n$. Between two candidates say Weibull-$R_1$ and Weibull-$R_2$ distributions, with respective pdfs $f_{R_1}$, $f_{R_2}$ and respective cdfs $F_{R_1}$, $F_{R_2}$, we decide in favor of one of them on the basis of the difference $\mathscr{D}_{1,2} = \eta_{\text{Weibull}-R_1} - \eta_{\text{Weibull}-R_2}$, where $\eta_{\text{Weibull}-R_1}$ and $\eta_{\text{Weibull}-R_2}$ are the entropies respectively of Weibull $- R_1$ and Weibull $- R_2$ random variable. We see that

$$\begin{aligned}\mathscr{D}_{1,2} &= (c-1)E\left[\log\frac{\log(1-F_{R_2}(X))}{\log(1-F_{R_1}(X))}\right] + E\left[\log\frac{f_{R_2}(X)}{1-F_{R_2}(X)}\frac{1-F_{R_1}(X)}{f_{R_1}(X)}\right] \\ &+ E\left\{\left[\frac{-\log(1-F_{R_1}(X))}{\gamma}\right]^c\right\} - E\left\{\left[\frac{-\log(1-F_{R_2}(X))}{\gamma}\right]^c\right\}.\end{aligned}$$

Large (respectively, small) values of $\mathscr{D}_{1,2}$ will support the Weibull-$R_1$ (respectively, Weibull-$R_2$) distribution.

**Proposition 4.** For fixed $c > 1$, $\gamma > 0$ and Weibull-$R_1$ stochastically larger than Weibull-$R_2$, $\mathscr{D}_{1,2}$ will support Weibull-$R_2$.

*Proof.* Since Weibull-$R_1$ is stochastically larger than Weibull-$R_2$, we can write $1 - F_{R_1} \geq 1 - F_{R_2}$. This implies the following

- $\left[\frac{-\log(1-F_{R_1}(X))}{\gamma}\right]^c \leq \left[\frac{-\log(1-F_{R_2}(X))}{\gamma}\right]^c$,

- $\log\left[\frac{1-F_{R_2}(X)}{1-F_{R_1}(X)}\right] \leq 0$.

So, $\mathscr{D}_{1,2}$ will be small which implies the result. $\square$

## 4.3 Moments

For any $r \in \mathbb{N}$, we can express the $r$th moment of $X$ as

$$\begin{aligned}E(X^r) &= \int_{-\infty}^{\infty} \frac{c}{\gamma} x^r \frac{f_R(x)}{1-F_R(x)}\left[\frac{-\log(1-F_R(x))}{\gamma}\right]^{c-1} \exp\left\{-\left[\frac{-\log(1-F_R(x))}{\gamma}\right]^c\right\} dx \\ &= \int_0^{\infty} e^{-u}\left[F_R^{-1}\left(1-\exp\left(-\gamma u^{1/c}\right)\right)\right]^r du,\end{aligned}$$

where $u = \left[\frac{-\log(1-F_R(x))}{\gamma}\right]^c$.

Analytical expressions for the moments do not appear possible, even for the four examples presented in Section 3. The moments should be computed numerically, for example, using integrate in the R software.

## 4.4 Reliability parameter

The reliability parameter $\mathcal{R}$ is defined as $\mathcal{R} = P(X > Y)$, where $X$ and $Y$ are independent random variables. Estimation of $\mathcal{R}$ is known as stress strength modeling. It has applications in many areas including break down of systems having two components. Other applications can be found in Weerahandi and Johnson (1992).

If $X$ and $Y$ are independent random variables with respective cdfs $F_1(x)$, $F_2(y)$ and respective pdfs $f_1(x)$, $f_2(y)$ then $\mathcal{R}$ can be written as

$$\mathcal{R} = P(X > Y) = \int_{-\infty}^{\infty} F_2(t) f_1(t) dt.$$

**Theorem 2.** Suppose $X$ and $Y$ are independent Weibull-$R\,(c_1, \gamma)$ and Weibull-$R\,(c_2, \gamma)$ random variables. Then

$$\mathcal{R} = 1 - \sum_{k=0}^{\infty} \frac{(-1)^k}{k!} \Gamma\left(k c_2 / c_1 + 1\right).$$

*Proof:* By (4) and (5),

$$\begin{aligned}
\mathcal{R} &= \int_{-\infty}^{\infty} \left[1 - \exp\left\{-\left[\frac{-\log(1 - F_R(x))}{\gamma}\right]^{c_2}\right\}\right] \\
&\quad \cdot \frac{c_1}{\gamma} \frac{f_R(x)}{1 - F_R(x)} \left[\frac{-\log(1 - F_R(x))}{\gamma}\right]^{c_1 - 1} \exp\left\{-\left[\frac{-\log(1 - F_R(x))}{\gamma}\right]^{c_1}\right\} dx \\
&= 1 - \int_0^{\infty} \exp(-u) \exp\left(-u^{c_2/c_1}\right) du \\
&= 1 - \sum_{k=0}^{\infty} \frac{(-1)^k}{k!} \Gamma\left(k c_2 / c_1 + 1\right),
\end{aligned}$$

where $u = \left[\frac{-\log(1 - F_R(x))}{\gamma}\right]^{c_1}$. Hence, the proof. $\square$

## 5 Conclusions

In this paper, we have introduced the Weibull-$R$ family with a hope that it will have more flexibility in situations where Weibull and other Weibull mixture distributions do not provide satisfactory fits. For each baseline distribution of $R$, our results can be easily adapted to obtain main structural properties of the Weibull-$R$ distribution. We have derived various properties of the Weibull-$R$ distributions, including the reliability parameter and the $r$th generalized moment. The proposed family unifies several previously proposed families of distributions, therefore yielding a general overview of these families for theoretical studies. It also provides a rather flexible mechanism for fitting a wide spectrum of real world data sets. For example, a Weibull-$R$ mixture distribution may be useful in the following scenarios:

- To characterize end-to-end Internet delay at coarse time-scales (Hernandez et al. (2006)).

- It provides a suitable distributions for modeling dependent lifetimes from heterogenous populations, as mixtures of defective devices with shorter lifetimes and standard devices with longer lifetimes.

- When $R$ the baseline distribution is Gompertz, a mixture of Weibull-Gompertz (in particular, the survival function) distribution will represent a theoretically motivated model for the scenario in which death or cases of a specific disease in an actual population can be due to sufficient causes from group 1 or group 2. For details, see Levy et al. (2014).

We hope that this family may attract wider applications in reliability and biology.

# References

[1] Ahsanullah, M. (1995). *Record Statistics.* Nova Science Publishers, Commack, New York.

[2] Alzaghal, A., Ghosh, I., and Alzaatreh, A. (2016). On Shifted Weibull-Pareto Distribution. *International Journal of Statistics and Probability,* 5, 139-149.

[3] Alzaatreh, A., and Ghosh, I. (2015). On the Weibull-X family of distributions. *Journal of Statistical Theory and Applications,* 14, 169-183.

[4] Alzaatreh, A., Famoye, F. and Lee, C. (2013a). Weibull-Pareto distribution and its applications. *Communications in Statistics - Theory and Methods*, 42, 1673-1691.

[5] Alzaatreh, A., Lee, C. and Famoye, F. (2013b). A new method for generating families of continuous distribution. *Metron*, 71, 63-79.

[6] Bain, L.J. (1978). *Statistical Analysis of Reliability and Life Testing Models.* Marcel Dekker, Inc, New York.

[7] Glick, N. (1978). Breaking record and breaking boards. *American Mathematical Monthly*, 85, 2-26.

[8] Gupta, R.D. and Kundu, D. (1997). Exponentiated exponential family: An alternative to gamma and Weibull distribution. *Technical Report*, Department of Mathematics, Statistics and Computer Science, University of New Brunswick, Saint-John, NB, Canada.

[9] Gurvich, M.R., Dibenedetto, A.T. and Rande, S.V. (1997). A new statistical distribution for characterizing the random strength of brittle materials. *Journal of Materials Science*, 32, 2559-2564.

[10] Hernandez, J.A. and Phillips, I.W. (2006). Weibull mixture model to characterise end-to-end Internet delay at coarse time-scales. *IEE Proceedings - Communications,* 153, 295-304.

[11] Huillet, T. and Raynaud, H.F. (1999). Rare events in a log-Weibull scenario: Application to earthquake magnitude data. *European Physics Journal B*, 12, 457-469.

[12] Lai, C.D., Murthy, D.N.P. and Xie, M. (2011). Weibull distributions. *Wiley Interdisciplinary Reviews*, 3, 282-287.

[13] Lee, C., Famoye, F., and Olumolade, O. (2007). Beta-Weibull distribution: Some properties and applications to censored data. *Journal of Modern Applied Statistical Methods*, 6, 173-186.

[14] Levy, G. and Levin, B. (2014). *The Biostatistics of aging. From Gompertzian Mortality to an Index of Aging-Relatedness.* John Wiley, NY.

[15] Mudholkar, G.S. and Kollia, G.D. (1994). Generalized Weibull family: A structural analysis. *Communications in Statistics - Theory and Methods*, 23, 1149-1171.

[16] Mudholkar, G.S. and Srivastava, D.K. (1993). Exponentiated Weibull family for analyzing bathtub failure-rate data. *IEEE Transactions on Reliability*, 42, 299-302.